\documentclass[11pt, reqno]{amsart}
\usepackage{amsmath}
\usepackage{amssymb}
\usepackage{amsthm}
\usepackage{enumerate}
\usepackage{xcolor}
\usepackage{mathrsfs}
\usepackage{hyperref}
\usepackage[abbrev]{amsrefs}
\usepackage[utf8]{inputenc}
\usepackage{ifthen}
\usepackage{enumitem}
\usepackage{bm}
\usepackage{graphicx}
\usepackage[all]{xy}

\usepackage{mathrsfs}
\newtheorem{thm}{}[section]
\newtheorem{theorem}[thm]{Theorem}

\newtheorem{lemma}[thm]{Lemma}

\theoremstyle{definition}

\theoremstyle{remark}
\newtheorem{remark}[thm]{Remark}

\numberwithin{equation}{section}
\allowdisplaybreaks

\newcommand{\N}{\mathbb{N}}

\newcommand{\R}{\mathbb{R}}
\newcommand{\C}{\mathbb{C}}

\renewcommand{\epsilon}{\varepsilon}

\newcommand{\DD}{\mathbb{D}}

\AtBeginDocument{
\def\MR#1{}
}

\title[Mean convergence]{Mean convergence of Fourier-Akhiezer-Chebyshev series}

\author[M. Bello-Hernández]{Manuel Bello-Hern\'{a}ndez}
\address{
Departamento de Matem\'aticas y Computaci\'on\\
Universidad de La Rioja\\
c/Madre de Dios, 53, Spain}
\email{mbello@unirioja.es}

\author[A. del Campo-López]{Alejandro del Campo López}
\address{
Universidad de La Rioja\\
c/Madre de Dios, 53, Spain}
\email{alejandro.del-campo@alum.unirioja.es}

\date{\today}

\begin{document}
     

	\begin{abstract}
	We prove mean convergence of the Fourier series in Akhiezer-Chebyshev polynomials in $L^p$, $p>1$, using a weighted inequality for the Hilbert transform in an arc of the unit circle.
	\end{abstract}
	
		\maketitle
	
\section{Introduction}

One of the cornerstone results in harmonic analysis has been the continuity of the conjugate operator or Hilbert transform in $L^p(0,2\pi)$ for $1<p<\infty$ by Marcel Riesz \cite{Rie27}. This result implies the convergence of the Fourier series in the mentioned function space. Riesz's results opened several questions of research such as finding closed systems of  functions in $L^p$.

The goal of this paper is to prove mean convergence of Fourier-Akhiezer-Chebyshev series. The Akhiezer-Chebyshev polynomials were introduced by  Akhiezer \cite{Akh60} and a rigorous exposition was given by Golinskii \cite{Gol98}. These polynomials are orthogonal with respect to a measure supported on an arc of the unit circle. At the present time nothing is known about convergence of
orthogonal Fourier series in $L^p$ spaces with arbitrary weights (see, for example, \cites{AskWai65,LeKusv99,MatNevTot1986,Muc1969,Muc70,Muc72,Pol47,Pol48}). Our results are limited to a type of Akhiezer-Chebyshev weights.

Set $\alpha\in(0,\pi),$ $\Delta_\alpha=\{e^{i\vartheta}:\alpha<\vartheta<2\pi-\alpha\}$.
Let $w_\alpha(\vartheta)$ be the Akhiezer-Chebyshev weight in $\Delta_\alpha$ given by
\begin{equation}
\label{eqAkhiezerWeight}
w_\alpha(\vartheta)=
\left\{
\begin{array}{ll}
\frac{\sin(\alpha/2)}{2\sin(\vartheta/2)\sqrt{\cos^2(\alpha/2)-\cos^2(\vartheta/2)}},& \vartheta\in(\alpha,2\pi-\alpha),\\
0,& \text{in other case}.
\end{array}
\right.
\end{equation}

Let $\{\psi_{n}\}_{n=0}^\infty$ be the sequence of orthonormal polynomials with
the inner product defined by 
\[
\langle f,g \rangle=\frac{1}{2\pi}\int_\alpha^{2\pi-\alpha}f(e^{i\vartheta})\overline{g(e^{i\vartheta})}w_\alpha(\vartheta)\, d\vartheta,
\]
Consider the Lebesgue space $L^p(w_\alpha)=L^p((\alpha,2\pi-\alpha),w_\alpha)$, $1\le p<\infty$, i.e. the set of all Borel measurable functions $f:\Delta_\alpha\to\C$ such that
\[
\int_\alpha^{2\pi-\alpha}|f(e^{i\vartheta})|^pw_\alpha(\vartheta)\, d\vartheta<\infty
\]
For $f\in L^p(w_\alpha)$, consider its $n-$th partial Fourier sum in terms of $\{\psi_{n}\}_{n=0}^\infty$
\[
\mathscr{S}_n(f,z)= \sum_{j=0}^n\langle f,\psi_j\rangle \psi_{j}(z).
\]
Then we prove the following result.

\begin{theorem}\label{TeoConvMedSerieFourierPsi} Let $p\in(1,\infty)$, and $f\in L^p(w_\alpha)$. Then we have
\begin{equation}
\label{eqTeoConvMedSerieFourierPsi}
\lim_{n\to\infty}\int_{\alpha}^{2\pi-\alpha}|f(e^{i\vartheta})-\mathscr{S}_n(f,e^{i\vartheta})|^pw_\alpha(\vartheta)\,d\vartheta=0.
\end{equation}

\end{theorem}

In our main auxiliary results we obtain a weighted inequality for
the Hilbert transform in an arc of the unit circle and we give an estimate for the sequence of para-orthonormal polynomials associated to Akhiezer-Chebyshev weight.

The paper is organized as follows. First, for an easy reading, we include some calculations for Akhiezer-Chebyshev polynomials in Section~\ref{sec:ACpoly}. Then, in Section~\ref{proofTeo1} we prove the weighted inequality for the Hilbert transform and Theorem~\ref{TeoConvMedSerieFourierPsi}. Finally, Section \ref{SectionAdd} contains some results on mean convergence of series in orthogonal polynomial with respect to Akhiezer-Chebyshev weight multiplied by a function with nice properties.

\section{Akhiezer-Chebyshev polynomials}\label{sec:ACpoly}
Several results in this section can be founded in \cite{Gol98}. We include the proofs of the main results in order to make the paper self-contained.
First, we consider a conformal mapping from the unit disk onto the complement of $\Delta_\alpha$. This function is written in terms of the product of two  M\"{o}bius transforms.

\begin{lemma}[\cite{Gol98}]\label{LemConfMap} Let $\beta=i\tan\frac{\pi-\alpha	}{4}$ and $V(z)=(z-e^{i\alpha})(z-e^{-i\alpha})$.

\begin{enumerate}
\item[(i)] The function 
\begin{equation}
\label{funcH}
z=h(v)=\frac{(v-\beta)(\beta v-1)}{(v+\beta)(\beta v+1)}
\end{equation}
is a conformal mapping from the unit disk $\mathbb{D}=\{z\in\C:|z|<1\}$ onto $\overline{\C}\setminus \Delta_{\alpha}$. Moreover, $e^{i\vartheta}=h(e^{i\omega})$ for $\omega\in[0,2\pi]$ sweeps the arc $e^{i\vartheta}\in\Delta_\alpha$ twice.

\item[(ii)] The functions
\begin{equation}
\label{eqV12}
\chi_{1,2}(z)=\frac{z+1\pm\sqrt{V(z)}}{2\cos(\alpha/2)}
\end{equation}
are conformal mappings from $\overline{\C}\setminus \Delta_\alpha$ to $\mathbb{E}=\{v\in\C:|v|>1\}$ or $\mathbb{D}$ according to either sign plus or minus is chosen, where 
we take the branch of the root such that $\sqrt{1}=1$.
\item[(iii)] Let $f$ be a Borel function.  We have $f\in L^1(w_\alpha)$ if and only if $f\circ h\in L^1(0,\pi)$. In that case, we get
\[
\int_\alpha^{2\pi-\alpha}f(e^{i\vartheta})w_\alpha(\vartheta)\,d\vartheta =\int_0^{\pi}f(h(e^{i\omega}))d\omega.
\]
\end{enumerate}
\end{lemma}
\begin{proof}
The M\"{o}bius transform
\begin{equation}
\label{Mobius}
w=w(v)=i\frac{1-\beta v}{v+\beta}
\end{equation}
maps the exterior of the unit circle $\mathbb{E}$ onto the interior of the unit circle, whereas $w(\frac1{v})=i\frac{v-\beta }{\beta v+1}$ maps the interior of the unit circle $\mathbb{D}$ onto itself. The inverse transform of $w=w(v)$ is
\begin{equation}
\label{InvMobVvsW}
v=v(w)=\frac{i-\beta w}{w+i\beta},
\end{equation}
which maps $\mathbb{E}$ onto $\mathbb{D}$.

Observe that
\[
z=h(v)=\frac{(v-\beta)(\beta v-1)}{(v+\beta)(\beta v+1)}=w(v)\,w(\frac1v).
\]
So we have
\begin{equation}
\label{symProp}
h(\frac1v)=h(v),\quad h(\frac1{\overline{v}})=\frac1{\overline{h(v)}},\quad h(\overline{v})=\frac{1}{\overline{h(v)}}.
\end{equation}
In particular,
\begin{equation}
\label{symPropCirc}
h(e^{-i\omega})=h(e^{i\omega}), \quad |h(e^{i\omega})|=1,\quad \omega\in\R.
\end{equation}
As both numerator and denominator of $h$ are polynomials of degree $2$, each point of $z\in\overline{\C}$ is the image by $z=h(v)$ of exactly two points. Since $h(1/v)=h(v)$, $h$ is injective in the unit disk and conformal since the derivative of $h$ is
\[
h'(v)=2(\beta+\beta^{-1})\frac{v^2-1}{(v+\beta)^2(v+\beta^{-1})^2}\ne 0, \quad v\in\mathbb{D}.
\]
Moreover, 
\[
h(1)=-\left(\frac{1-\beta}{1+\beta}\right)^2=e^{i\alpha},\quad h(-1)=e^{-i\alpha},\quad h(\pm i)=-1,\quad h(0)=h(\infty)=1,
\]
\[
|h(v)|=1,\quad v\in\R.
\]
Since
\[
h(i\epsilon)=\frac{1-2\tan(\alpha/2)\epsilon-\epsilon^2}{1+2\tan(\alpha/2)\epsilon-\epsilon^2},
\]
for $\epsilon>0$ small, we have $|h(i\epsilon)|<1$. Hence,
$z=h(v)$ maps $\DD$ onto $\overline{\C}\setminus \Delta_\alpha$; and $e^{i\vartheta}=h(e^{i\omega})$ for $\omega\in[0,2\pi]$ sweeps the arc $\Delta_\alpha$ twice.

Next we get (\textit{ii}). Taking into account \eqref{InvMobVvsW}, the composition $z=h(v(w))$ maps $\mathbb{E}$ onto $\overline{\C}\setminus\Delta_\alpha$ given by
\begin{equation}
\label{eqHmapW}
z=h(v(w))=\frac{\cos(\alpha/2)w^2-w}{w-\cos(\alpha/2)},
\end{equation}
which has the property
\[
\lim_{w\to\infty}\frac{h(v(w))}{w}=\cos(\alpha/2).
\]
Then the functions given in \eqref{eqV12}, which are the solutions of \eqref{eqHmapW} for each $z$, are conformal  mappings from either $\mathbb{E}$ or $\mathbb{D}$ to $\overline{\C}\setminus \Delta_\alpha$ according to whether $\sqrt{1}=1$ or $\sqrt{1}=-1$, respectively.
The number $\gamma=\cos(\alpha/2)$ is the transfinite diameter of $\Delta_\alpha$. Observe that
\[
V(z)=(z+1)^2-4z\cos^2(\alpha/2)=z^2-2z\cos\alpha+1.
\]

Now we check (\emph{iii}). If $e^{i\vartheta}=h(e^{i\omega})$, then 
\begin{equation}
\label{eqRelOmegaTheta}
e^{2i\omega}-i(\beta+\beta^{-1})\cot \frac\vartheta{2}e^{i\omega}+1=0,
\end{equation}
and this equation for $\omega\in(0,\pi)$ (i.e. $\sin \omega>0$) has solution
\begin{equation}
\label{eqOmegaTheta}
e^{i\omega}=\frac{\tan\frac\alpha2}{\tan\frac\vartheta2} +i\sqrt{1-\left(\frac{\tan\frac\alpha2}{\tan\frac\vartheta2}\right)^2}
\end{equation}
and for $\omega\in(\pi,2\pi)$,
\[
e^{i\omega}=\frac{\tan\frac\alpha2}{\tan\frac\vartheta2} -i\sqrt{1-\left(\frac{\tan\frac\alpha2}{\tan\frac\vartheta2}\right)^2}.
\]
Thus, we have
\[
2i e^{2i\omega}\frac{d\omega}{d\vartheta}-i(\beta+\beta^{-1})(-\frac1{2\sin^2(\vartheta/2)}e^{i\omega} +i\cot(\vartheta/2)e^{i\omega}\frac{d\omega}{d\vartheta})=0,
\]
and
\begin{equation}
\label{eqDerOmegaTheta}
\frac{d\omega}{d\vartheta}= \frac{\sin(\alpha/2)}{2\sin(\vartheta/2)\sqrt{\cos^2(\alpha/2)-\cos^2(\vartheta/2)}},\quad \vartheta\in (\alpha,2\pi-\alpha).
\end{equation} 
Therefore, statement (\emph{iii}) follows straightforward.
\end{proof}

%

Let $K=\sqrt{\frac{2\sin(\alpha/2)}{1+\sin(\alpha/2)}}$. The next lemma gives the most important relations of such polynomials for our interest in this paper.
\begin{lemma}[\cite{Gol98}]\label{lemAux}  For $z\in \C\setminus\Delta_\alpha$ we have
\begin{equation}\label{eqAkhChebPolys1}
K^{-1}\psi_n(z)=\frac{i\beta}{1+\beta^2}(\chi_1(z)^{n-1}+\chi_2(z)^{n-1})+\frac{1}{1+\beta^2}(\chi_1(z)^{n}+\chi_2(z)^{n}),\quad n\ge 1.
\end{equation}
Also, for $\vartheta\in(\alpha,2\pi-\alpha)$ we have
\begin{equation}\label{eqAkhChebPolys}
K^{-1}\psi_n(e^{i\vartheta})=2e^{in\vartheta/2}\left\{\frac{i\beta}{1+\beta^2}e^{-i\vartheta/2}\cos((n-1)\lambda)+\frac{1}{1+\beta^2}\cos(n\lambda)\right\},
\end{equation}
\begin{equation}
\label{eqPolRev}
K^{-1}\psi_n^*(e^{i\vartheta})=2e^{in\vartheta/2} \left\{\frac{i\beta}{1+\beta^2}e^{i\vartheta/2}\cos((n-1)\lambda)+\frac{1}{1+\beta^2}\cos(n\lambda)\right\}, \quad n\ge 1,
\end{equation}
where $
\cos\lambda=\frac{\cos(\vartheta/2)}{\cos(\alpha/2)},\,\lambda\in[0,\pi], \, \vartheta\in[\alpha,2\pi-\alpha].$
Moreover, there exists a constant $C>0$ such that\footnote{We use the same letter $C$ for different constants independent of either the counter $n$, a point in an interval or the function in a specified class. Even the appearance of  the same letter $C$ in two consecutive inequalities  can refer to different values.}
\begin{equation}
\label{eqBoundAkh}
|\psi_n(e^{i\vartheta})|\le C,\quad \text{for all }\vartheta\in [\alpha,2\pi-\alpha].
\end{equation}
\end{lemma}

\begin{proof} We have 
\[
\chi_1(z)+\chi_2(z)=\frac{z+1}{\cos(\alpha/2)},\quad \chi_1(z)\chi_2(z)=z,
\]
Thus, for each $z\in \overline{\C}\setminus \Delta_\alpha$ there exists $v$ with $|v|<1$ such that
\[
\chi_1(z)=w(v),\quad \chi_2(z)=w(1/v).
\]
Let $g_n(z)$ be the hand right side of \eqref{eqAkhChebPolys1}. 
 By \eqref{eqV12}, we get
\begin{multline*}
\chi_1^{n+1}(z)+\chi_2^{n+1}(z)
\\= \left(\chi_1^{n}(z)+\chi_2^{n}(z)\right)\frac{z+1}{\cos(\alpha/2)}-z \left(\chi_1^{n-1}(z)+\chi_2^{n-1}(z)\right),
\end{multline*}
which by induction yields that $g_n(z)$ is a polynomial of degree $n$ in $z$.
Moreover, 
\begin{equation}
\label{eqCoefCons}
\frac{1}{1-\beta v}=i \frac{\beta}{1+\beta^2}\frac{1}{w(v)}+\frac{1}{1+\beta^2},
\quad
\frac{v}{v-\beta}=i \frac{\beta}{1+\beta^2}\frac{1}{w(1/v)}+\frac{1}{1+\beta^2}.
\end{equation}
Thus,
\[
g_n(z)=\frac{1}{1-\beta v}w(v)^n+\frac{v}{v-\beta}w(1/v)^n.
\]
Since $z=h(v)=w(v)w(1/v)$, the leading coefficient $\alpha_n$ of $g_n(z)$ is
\begin{multline*}
\alpha_n=\lim_{z\to\infty}\frac{K g_n(z)}{z^n}=\lim_{v\to-\beta}\left(\frac{v}{v-\beta}\frac{1}{w(v)^n}+\frac{1}{1-\beta v}\frac{1}{w(1/v)^n}\right)\\=K\frac{1+\sin(\alpha/2)}{2\sin(\alpha/2)}\frac{1}{\cos^n(\alpha/2)}=\frac{1}{K\cos^n(\alpha/2)}.
\end{multline*}
From Lemma~\ref{LemConfMap}~(\emph{iii}) and Cauchy's theorem, for $m<n$ we obtain 
\begin{multline*}
\int_{\alpha}^{2\pi-\alpha}g_n(e^{i\vartheta})e^{-im\vartheta}w_\alpha(\vartheta)\,d\vartheta\\= \int_0^{\pi} \left\{\frac{w(e^{i\omega})^n}{1-\beta e^{i\omega}}+ \frac{e^{i\omega}w(e^{-i\omega})^n}{e^{i\omega}-\beta}\right\}\frac1{h(e^{i\omega})^m}\,d\omega\\= i^{n-2m}\int_0^{\pi} \left(\frac{1-\beta e^{i\omega}}{e^{i\omega}+\beta}\right)^{n-m}\left(\frac{1+\beta e^{i\omega}}{e^{i\beta}-\beta}\right)^m\frac{d\omega}{(1-\beta e^{i\omega})}\\ + i^{n-2m}\int_0^\pi\left(\frac{e^{i\omega}-\beta}{1+\beta e^{i\omega}}\right)^{n-m}\left(\frac{e^{i\omega}+\beta}{1-\beta e^{i\omega}}\right)^{m}\frac{e^{i\omega}d\omega}{(e^{i\omega}-\beta)}\\= i^{n-2m-1}\int_{\mathbb{T}}\left(\frac{\zeta-\beta}{1+\beta \zeta}\right)^{n-m}\left(\frac{\zeta+\beta}{1-\beta \zeta}\right)^{m}\frac{d\zeta}{(\zeta-\beta)}=0.
\end{multline*}
The same calculation for $m=n$ yields
\begin{multline*}
\int_{\alpha}^{2\pi-\alpha}g_n(e^{i\vartheta})e^{-in\vartheta}w_\alpha(\vartheta)d\vartheta=i^{-n-1}\int_{\mathbb{T}}\left(\frac{\zeta+\beta}{1-\beta \zeta}\right)^{n}\frac{d\zeta}{(\zeta-\beta)}\\=i^{-n}2\pi  \left(\frac{2\beta}{1-\beta^2}\right)^{n}.
\end{multline*}
Therefore, with the correct normalization we obtain \eqref{eqAkhChebPolys1}.

Moreover, since $V(0)=1$ and
\begin{multline*}
\chi_{1,2}(e^{i\vartheta})=e^{i\vartheta/2}\frac{\cos(\vartheta/2)}{\cos{\alpha/2}}\pm\frac{1}{\cos(\alpha/2)} \sqrt{e^{i\vartheta}i^2\sin(\frac{\vartheta-\alpha}{2})\sin(\frac{\vartheta+\alpha}{2})}\\= e^{i\vartheta/2}\frac{\cos(\vartheta/2)}{\cos(\alpha/2)}\pm e^{i\vartheta/2}i\sqrt{(1-\frac{\cos^2(\vartheta/2)}{\cos^2(\alpha/2)}})=e^{i\vartheta/2\pm\lambda i},
\end{multline*}
where $\cos\lambda=\frac{\cos(\vartheta/2)}{\cos(\alpha/2)}$ for $\vartheta\in(\alpha,2\pi-\alpha)$,  we get
\[
\chi_1(z)^j+\chi_2(z)^j=2e^{ij\vartheta/2}\cos(j\lambda ),\quad j\ge 1.
\]
 Therefore, \eqref{eqAkhChebPolys} follows intermediately from \eqref{eqAkhChebPolys1}.

Statement \eqref{eqBoundAkh} is clear from \eqref{eqAkhChebPolys}. The relation in \eqref{eqPolRev} follows from the definition of the reverse polynomial.
\end{proof}

Next we consider para-orthogonal polynomials associated to Akhiezer-Chebyshev weight. These polynomials allow us to work on the diagonal of integration region to estimate the $L^p(w_\alpha)-$norm of $\mathscr{S}_n$.
The para-orthogonal polynomial associated to Akhiezer-Chebyshev polynomials with a zero at $e^{\pm i\alpha}$ is defined by 
\[
\Lambda_{n}^{(\pm\alpha)}(z)=\psi_n^*(z)-\frac{\psi_n^*(e^{\pm i\alpha})}{\psi_n(e^{\pm i\alpha})}\psi_n(z).
\]
Observe that $\Lambda_{n}^{(\pm\alpha)}(e^{\pm i\alpha})=0$ and $|\frac{\psi_n^*(e^{\pm i\alpha})}{\psi_n(e^{\pm i\alpha})}|=1$, so by \eqref{eqBoundAkh} the sequence of para-orthogonal polynomials $\{\Lambda_{n}^{(\pm\alpha)}\}_{n=0}^{\infty}$ is uniformly bounded in the $\Delta_\alpha$.  Moreover, we can say much more with the following estimate.
\begin{lemma}\label{LemaAcotLambda} There exists a constant $C$ independent of $n$ and $\vartheta$ such that
\begin{equation}
\label{estParaOrth}
|\Lambda_{n}^{(\pm\alpha)}(e^{i\vartheta})|\le C|e^{i\vartheta}-e^{\pm i\alpha}|^{1/2},\quad \text{for all } n\in\N,\quad \vartheta\in(\alpha,2\pi-\alpha).
\end{equation}
\end{lemma}

\begin{proof} We give only the computation to show that there exists a constant $C>0$ such that
\begin{equation}
\label{estParaOrth+}
|\Lambda_{n}^{(\alpha)}(e^{i\vartheta})|\le C|e^{i\vartheta}-e^{i\alpha}|^{1/2},\quad \text{for all }n\in\N,\quad \vartheta\in(\alpha,2\pi-\alpha),
\end{equation}
since the other inequality is proved in the same manner.
Let $\lambda_\vartheta$ be the value in $[0,\pi]$ such that
\[
\cos\lambda_\vartheta=\frac{\cos(\vartheta/2)}{\cos(\alpha/2)},\quad \vartheta\in [\alpha,2\pi-\alpha].
\]
Of course, $\lambda_\alpha=0$, so 
\begin{equation}
\label{eqLambdan}
\cos(n\lambda_\alpha)=1,\quad \text{for all }n\in\N.
\end{equation} 
Moreover, from \eqref{eqAkhChebPolys},  \eqref{eqPolRev}, and \eqref{eqLambdan} it follows that
\begin{equation}
\label{Delta}
\frac{\psi_n^*(e^{i\alpha})}{\psi_n(e^{i\alpha})}=\frac{i\beta e^{i\alpha/2}+1}{i\beta e^{-i\alpha/2}+1}.
\end{equation}
Denote by $\Upsilon=\frac{\psi_n^*(e^{i\alpha})}{\psi_n(e^{i\alpha})}$ the value above. Observe that
\begin{equation}
\label{1mendosDelta}
1-\Upsilon=i\beta\frac{e^{-i\alpha/2}-e^{i\alpha/2}}{i\beta e^{-i\alpha/2}+1}.
\end{equation}
Thus from \eqref{eqAkhChebPolys} and \eqref{eqPolRev} we obtain
\begin{multline}\label{eqAux1}
\frac{K^{-1}}{2}|\Lambda_n^{(\alpha)}(e^{i\vartheta})|=\left|\frac{i\beta}{1+\beta^2} \cos((n-1)\lambda_\vartheta) (e^{-i\vartheta/2}-\Upsilon e^{i\vartheta/2})\right.\\\left.+\frac{1}{1+\beta^2}\cos(n\lambda_\vartheta) (1-\Upsilon)\right|\\ 
=\left|\cos((n-1)\lambda_\vartheta)\left( \frac{i\beta}{1+\beta^2}  (e^{-i\vartheta/2}-\Upsilon e^{i\vartheta/2})+\frac{\cos\lambda_\vartheta}{1+\beta^2}(1-\Upsilon)\right)\right.\\-\left.\frac{\sin \lambda_\vartheta\sin((n-1)\lambda_\vartheta)}{1+\beta^2}(1-\Upsilon)\right|.
\end{multline}
Since $
\sin \lambda_\vartheta=\frac{\sqrt{\cos^2(\alpha/2)-\cos^2(\vartheta/2)}}{\cos(\alpha/2)}$ and $ \frac{i\beta}{1+\beta^2}  (e^{-i\vartheta/2}-\Upsilon e^{i\vartheta/2})+\frac{\cos\lambda_\vartheta}{1+\beta^2}(1-\Upsilon)=e^{-i\vartheta/2}P(e^{i\vartheta})$, where $P$ is a polynomial of degree 1 which equals zero at $\vartheta=\alpha$, there exists a constant $C$ independent of $\vartheta$ such that
\[
\sin \lambda_\vartheta\le C|e^{i\vartheta}-e^{i\alpha}|^{1/2},
\]
and
\[
\left| \frac{i\beta}{1+\beta^2}  (e^{-i\vartheta/2}-\Upsilon e^{i\vartheta/2})+\frac{\cos\lambda_\vartheta}{1+\beta^2}(1-\Upsilon)\right|\le C|e^{i\vartheta}-e^{i\alpha}|^{1/2}.
\]
Therefore, \eqref{estParaOrth+} follows immediately from \eqref{eqAux1} and the inequalities above.
\end{proof}

\section{Proof of Theorem \ref{TeoConvMedSerieFourierPsi}}\label{proofTeo1}

 We need several results on the Hilbert transform on an arc of the unit circle. Let us recall Riesz's theorem \cite{Rie27} on this issue. This result states that the harmonic conjugate
\[
\overset{\backsim}{f}(x)=\frac{1}{2\pi}\textnormal{PV}\int_0^{2\pi} f(t)\cot\frac{x-t}{2}\, dt,\quad x\in (0,2\pi),
\]
is a bounded operator from $L^p(0,2\pi)$ to itself, where the integral is defined as the Cauchy principal value at $t=x$; it means that $\overset{\backsim}{f}$ exists almost everywhere and there exists a constant $C>0$ such that for all $f\in L^p(0,2\pi)$, $p>1$, it holds
\begin{equation}
\label{bouHilbTrans}
\int_{0}^{2\pi}|\overset{\backsim}{f}(x)|^p\, dx\le C \int_{0}^{2\pi}|f(x)|^p\, dx.
\end{equation}
Observe that if $f$ is a $2\pi-$periodic even function, then 
\[
\overset{\backsim}{f}(x)=\frac{1}{2\pi}\textnormal{PV}\int_0^{\pi} f(t)(\cot\frac{x+t}{2}-\cot\frac{t-x}{2})\, dt.
\]
In this case, Riesz's result \eqref{bouHilbTrans} can be written as
\begin{equation}
\label{bouHilbTrans2}
\int_{0}^{\pi}|\overset{\backsim}{f}(x)|^p\, dx\le C \int_{0}^{\pi}|f(x)|^p\, dx,
\end{equation}
for all $f\in L^p(0,\pi)$.

Let $f\in L^p(w_\alpha)$. The (unweighted) Hilbert transform in the arc $\Delta_\alpha$ is defined by
 \begin{multline*}
\mathcal{H}_1(f)(e^{i\tau})=\frac{1}{2\pi}\textnormal{PV}\int_\alpha^{2\pi-\alpha}\frac{f(e^{i\vartheta})}{e^{i\tau}-e^{i\vartheta}}\, d\vartheta\\
=\frac{1}{2\pi}\lim_ {\epsilon\to 0^+}\int_{(\alpha,2\pi-\alpha)\setminus(\tau-\epsilon,\tau+\epsilon)}\frac{f(e^{i\vartheta})}{e^{i\tau}-e^{i\vartheta}}\, d\vartheta,
\end{multline*}
where $\tau\in(\alpha,2\pi-\alpha)$. The topic of Hilbert transform along curves has a long history in harmonic
analysis that it can be seen in \cite{SteWai78}. To get an weighted inequality for this transform we shall use the following Muckenhoupt's inequality (see \cite[Lemma 8, p. 440 of the second paper]{Muc70}).
\begin{lemma}
If $1<p<\infty$, $r>-1/p$, $s<1-1/p$, $R<1-1/p$, $S>-1/p$, $r\ge R$ and $s\le S$, then there exists a constant $C$, independent of $f$, such that
\begin{equation}
\label{desMuc}
\int_{-\infty}^{\infty}\left|\int_{-\infty}^{\infty}\frac{f(y)}{x-y}|x|^r(1+|x|)^{s-r}\,dy\right|^p\,dx\le C\int_{-\infty}^{\infty}\left|f(y)|y|^R(1+|y|)^{S-R}\right|^p\, dy.
\end{equation}
\end{lemma}
Next we state  an  weighted inequality for the Hilbert transform in an arc. 

\begin{lemma}\label{MucOnAnArc} If $1<p<\infty$, then there exists a constant $C$ such that for all $f\in L^p(w_\alpha)$ we have
\begin{equation}
\label{IneqNonWeightedHilbertTransform}
\int_\alpha^{2\pi-\alpha}\left|\mathcal{H}_1(f)(e^{i\tau})\right|^pw_\alpha(\tau)\,d\tau\le C\int_\alpha^{2\pi-\alpha}|f(e^{i\tau})|^pw_\alpha(\tau)\, d\tau.
\end{equation}
In particular, $\mathcal{H}_1(f)(e^{i\tau})$ exists a.e. for all $f\in L^p(w_\alpha)$.
\end{lemma}

\begin{proof}
Easy computations give  us
\[
h(e^{it})-h(e^{is})=\frac{2(\beta+\beta^{-1})(e^{it}-e^{is})(e^{i(s+t)}-1)}{(e^{it}+\beta)(e^{is}+\beta)(e^{it}+\beta^{-1})(e^{is}+\beta^{-1})},
\] 
\[
e^{it}-e^{is}=2ie^{i(s+t)/2}\sin\frac{t-s}{2},
\]
 \[
 e^{i(s+t)}-1=2ie^{i(s+t)/2}\sin\frac{t+s}{2},
\]
and
\begin{multline}
\label{MainTransfH}
h(e^{it})-h(e^{is})=\frac{8(\beta+\beta^{-1})e^{i(s+t)}}{(e^{it}+\beta)(e^{is}+\beta)(e^{it}+\beta^{-1})(e^{is}+\beta^{-1})}\frac{\sin s}{\cot\frac{t+s}{2}-\cot\frac{t-s}{2}}\\
=\frac{4(\beta+\beta^{-1})e^{i(s+t)}(\cos s-\cos t)}{(e^{it}+\beta)(e^{is}+\beta)(e^{it}+\beta^{-1})(e^{is}+\beta^{-1})}.
\end{multline}
So,  we can write
\begin{equation}
\label{CambH1}
\frac{f(h(e^{it}))}{h(e^{is})-h(e^{it})}=\frac{f(h(e^{it}))(e^{it}+\beta)(e^{is}+\beta)(e^{it}+\beta^{-1})(e^{is}+\beta^{-1})}{4(\beta+\beta^{-1})e^{i(s+t)}(\cos s-\cos t)}.
\end{equation}
Using \eqref{eqOmegaTheta}, we also have
\begin{equation}
\label{FormulaSin}
\sin t=\frac{\sqrt{\cos^{2}(\alpha/2)-\cos^2(\vartheta/2)}}{\cos(\alpha/2)\sin(\vartheta/2)}=\frac{\tan(\alpha/2)}{2\sin^2(\vartheta/2) w_{\alpha}(\vartheta)}.
\end{equation}
Hence, doing the  change of variables $e^{i\tau}=h(e^{is})$, $e^{i\vartheta}=h(e^{it})$ and $s,t\in (0,\pi)$, from Lemma \ref{LemConfMap} (\emph{iii}), we get
\begin{multline}
\label{eqH1To0pi}
\int_\alpha^{2\pi-\alpha}|\mathcal{H}_1(f)(e^{i\tau})|^pw_\alpha(\tau)\, d\tau\\=\int_\alpha^{2\pi-\alpha}\left|\int_\alpha^{2\pi-\alpha}\frac{f(e^{i\vartheta})}{e^{i\tau}-e^{i\vartheta}}\, d\vartheta\right|^pw_\alpha(\tau)\, d\tau\\ =\int_0^{\pi}\left|\int_0^{\pi}\frac{2f(h(e^{it}))\sin t\cot(\alpha/2)\sin^2(\vartheta/2)}{h(e^{is})-h(e^{it})}\, dt\right|^p\, ds.
\end{multline}

Since there exist positive constants $C_1,C_2$ such that
\[
C_1\le \sin (\vartheta/2)\le C_2,\quad \vartheta\in(\alpha,2\pi-\alpha),
\]
\[
C_1\le \left|\frac{(e^{it}+\beta)(e^{is}+\beta)(e^{it}+\beta^{-1})(e^{is}+\beta^{-1})}{e^{i(s+t)}}\right|\le C_2,\quad s,t\in (0,\pi),
\]
and by \eqref{CambH1}  and \eqref{eqH1To0pi} we obtain that there exists a constant $C>0$ for all $f\in L^p(w_\alpha)$ such that
\begin{equation}
\label{eqHilbertTransf1}
\int_\alpha^{2\pi-\alpha}|\mathcal{H}_1(f)(e^{i\tau})|^pw_\alpha(\tau)\, d\tau\\\le C \int_0^\pi \left|\int_0^{\pi}\frac{f(h(e^{it}))\sin t}{\cos t-\cos s}\, dt\right|^p\, ds.
\end{equation}
Next, we change the variables of integration to $\sqrt{u}=\tan t/2$, $\sqrt{v}=\tan s/2$ in the hand right side of \eqref{eqHilbertTransf1},
\begin{multline}\label{eqAux2}
 C_1\int_0^\pi \left|\int_0^{\pi}\frac{f(h(e^{it}))\sin t}{\cos t-\cos s}\, dt\right|^p\, ds\\=  \int_0^\infty \left|\int_0^{\infty}\frac{f(h(e^{it}))}{v-u}\frac{1+v}{1+u}\, \frac{du}{1+u}\right|^p\, \frac{dv}{\sqrt{v}(1+v)}\\ \le C_2\left(\int_0^\infty \left|\int_0^{\infty}f(h(e^{it}))/(1+u)^2\,du\right|^p\, \frac{dv}{\sqrt{v}(1+v)}\right.\\+\left. \int_0^\infty \left|\int_0^{\infty}\frac{f(h(e^{it}))/(1+u)}{v-u}\, du\right|^p\, \frac{dv}{\sqrt{v}(1+v)}\right),
\end{multline}
where $C_1,C_2$ are constants independent of $f$ and in the inequality above we have used that $\frac{1+v}{1+u}=\frac{v-u}{1+u}+1$. 
Moreover, $f\in L^p((\alpha,2\pi-\alpha),w_\alpha)$ if and only if $f(h(e^{2i\arctan\sqrt{v}}))\in L^p((0,\infty),\frac{1}{\sqrt{v}(1+v)})$, and 
\[
\int_\alpha^{2\pi-\alpha}|f(e^{i\tau})|^pw_{\alpha}(\tau)\, d\tau=\int_0^{\infty}\left|f(h(e^{2i\arctan\sqrt{v}}))\right|^p\frac{dv}{\sqrt{v}(1+v)}.
\]
Obviously, in that case we also have  $ f(h(e^{i2i\arctan\sqrt{v}}))/(1+v)\in L^p((0,\infty),\frac{1}{\sqrt{v}(1+v)})$ and by H\"{o}lder's inequality there exists $C$ independent of $f$ such that
\begin{equation}
\label{eqAux3}
\left|\int_0^{\infty}f(h(e^{2i\arctan\sqrt{u}}))/(1+u)^2\,du\right|\le C
\left(\int_0^{\infty}\left|f(h(e^{2i\arctan\sqrt{u}}))\right|^p\frac{du}{\sqrt{u}(1+u)}\right)^{1/p}.
\end{equation}
Furthermore, applying Muckenhoupt's inequality \eqref{desMuc} to the function which equals  $f(h(e^{2i\arctan\sqrt{u}}))/(1+u) $ for $u>0$ and zero for $u\le 0$ with $s=-3/(2p)$, $S=1-3/(2p)$, $r=R=-1/(2p)$, we obtain
\begin{multline}
\label{eqAux4}
\int_0^\infty \left|\int_0^{\infty}\frac{f(h(e^{2i\arctan\sqrt{u}}))/(1+u)}{v-u}\, du\right|^p\, \frac{dv}{\sqrt{v}(1+v)}\\ \le C
\int_0^{\infty}\left|f(h(e^{2i\arctan\sqrt{v}}))\right|^p\frac{dv}{\sqrt{v}(1+v)}.
\end{multline}
Combining relations \eqref{eqHilbertTransf1}, \eqref{eqAux2}, \eqref{eqAux3}, and \eqref{eqAux4}, we get \eqref{IneqNonWeightedHilbertTransform}.
\end{proof}

Now we can prove Theorem \ref{TeoConvMedSerieFourierPsi}.

\begin{proof}[Poof of Theorem \ref{TeoConvMedSerieFourierPsi}]
If we prove that for each $f\in L^p(w_\alpha)$ there exists a constant $C=C(f)$ such that
\begin{equation}
\label{AcotaFourierSeriesPsi}
\int_\alpha^{2\pi-\alpha}|\mathscr{S}_n(f,e^{i\tau})|^pw_\alpha(\tau)\, d\tau\le C\int_\alpha^{2\pi-\alpha}|f(e^{i\tau})|^pw_\alpha(\tau)\, d\tau,
\end{equation}
for all $n\in\N$, then by the Banach-Steinhaus theorem actually there exists a constant $C$ independent of $f\in L^p(w_\alpha)$ such that the inequality above holds.
Once we know that the constant $C$ in \eqref{AcotaFourierSeriesPsi} is independent of $f$, by Szeg\H{o}-Kolmogorov-Krein's theorem  (\cite[Addenda B]{Akh92}, \cite[chapter 1]{GerBook61}, \cite[chapter 3]{GreSze84}) the algebraic polynomials are dense in $L^p(w_\alpha)$ and therefore from \eqref{AcotaFourierSeriesPsi} the statement in Theorem \ref{TeoConvMedSerieFourierPsi} follows.

Moreover, it is sufficient to prove \eqref{AcotaFourierSeriesPsi} with the integral on the left taken over $(\alpha,\pi)$ and on $(\pi,2\pi-\alpha)$; the second one is obtained with analogous arguments to the first one. So, for $1<p<\infty$ and $f\in L^p(w_\alpha)$, we  check only that there exists a constant $C=C(f,p)$ such that  we have
\begin{equation}
\label{AcotaFourierSeriesPsi1}
\int_\alpha ^{\pi}|\mathscr{S}_n(f,e^{i\tau})|^pw_\alpha(\tau)\, d\tau\le C\int_\alpha^{2\pi-\alpha}|f(e^{i\tau})|^pw_\alpha(\tau)\, d\tau,
\end{equation}
for all $n\in\N$.

According to Christoffel-Darboux formula it holds
 \[
K_n(\vartheta,\tau)=\sum_{j=0}^n\overline{\psi_{j}(e^{i\vartheta})}\psi(e^{i\tau})=\frac{\overline{\psi_{n+1}^*(e^{i\vartheta})}\psi_{n+1}^*(e^{i\tau})-\overline{\psi_{n+1}(e^{i\vartheta})}\psi_{n+1}(e^{i\tau})}{1-\overline{e^{i\vartheta}}e^{i\tau}},
\]
and
\begin{equation}\label{ChrDarForPsi}
\mathscr{S}_n(f,e^{i\tau})=\frac{1}{2\pi}\int_{\alpha}^{2\pi-\alpha} f(e^{i\vartheta})\overline{K_{n}(\tau,\vartheta)}w_{\alpha}(\vartheta)\,d\vartheta.
\end{equation}
Since
\[
\frac{\psi_{n+1}^*(e^{i\alpha})}{\psi_{n+1}(e^{i\alpha})}=\frac{e^{in\alpha}\overline{\psi_{n+1}(e^{i\alpha})}}{\psi_{n+1}(e^{i\alpha})}=\frac{\overline{\psi_{n+1}(e^{i\alpha})}}{\overline{\psi_{n+1}^*(e^{i\alpha})}},
\]
we have
\begin{equation}
\label{eqDaigonalDecomposition}
(1-\overline{e^{i\vartheta}}e^{i\tau})K_{n}(\vartheta,\tau)=\overline{\psi_{n+1}^*(e^{i\vartheta})}\Lambda_{n+1}^{(\alpha)}(e^{i\tau})+\psi_{n+1}(e^{i\tau}) \frac{\psi_{n+1}^*(e^{i\alpha})}{\psi_{n+1}(e^{i\alpha})}\overline{\Lambda_{n+1}^{(\alpha)}(e^{i\vartheta})}.
\end{equation}
Let $\delta>0$ be small enough (smaller than $\pi-\alpha$). We have
\begin{multline}
\label{DescSn1}
\mathscr{S}_n(f,e^{i\tau})=
\frac{1}{2\pi}\int_{\alpha}^{\pi+\delta} f(e^{i\vartheta})K_{n}(\vartheta,\tau)w_\alpha(\vartheta)\,d\vartheta\\+\frac{1}{2\pi}\int_{\pi+\delta}^{2\pi-\alpha} f(e^{i\vartheta})K_{n}(\vartheta,\tau)w_\alpha(\vartheta)\,d\vartheta.
\end{multline}

The integral 
\[
\int_\alpha^{\pi}\left|\frac{1}{2\pi}\int_{\pi+\delta}^{2\pi-\alpha} f(e^{i\vartheta})K_n(\vartheta,\tau)w_\alpha(\vartheta)\,d\vartheta\right|^pw_\alpha(\tau)\, d\tau
\]
is not a singular integral and this is a bounded sequence because of the sequence of polynomials $\{\psi_n\}_{n=0}^{\infty}$ also is. Thus, there exists a constant $C>0$ such that
\begin{multline}
\label{ActCero}
\int_\alpha^{\pi}\left|\frac{1}{2\pi}\int_{\pi+\delta}^{2\pi-\alpha} f(e^{i\vartheta})K_n(\vartheta,\tau)w_\alpha(\vartheta)\,d\vartheta\right|^pw_\alpha(\tau)\,d\tau\\ \le C \int_\alpha^{2\pi-\alpha}|f(e^{i\tau})|^pw_\alpha(\tau)\, d\tau.
\end{multline}
Therefore, it is enough to work with
\[
\int_\alpha^{\pi}\left|\int_\alpha^{\pi+\delta} f(e^{i\vartheta})K_n(\vartheta,\tau)\,w_\alpha(\vartheta)\,d\vartheta\right|^pw_\alpha(\tau)\, d\tau.
\]

By \eqref{eqDaigonalDecomposition} we obtain
\begin{multline}
\label{RelaBasicaDiagonal}
\int_\alpha^{\pi}\left|\int_\alpha^{\pi+\delta} f(e^{i\vartheta})K_n(\vartheta,\tau)w_\alpha(\vartheta)\,d\vartheta\right|^pw_\alpha(\tau)\\\le C
\left(\int_\alpha^{\pi}\left|\Lambda_{n+1}^{(\alpha)}(e^{i\tau})\int_\alpha^{\pi+\delta} f(e^{i\vartheta})\frac{\overline{\psi_{n+1}^*(e^{i\vartheta)})}}{1-\overline{e^{i\vartheta}}e^{i\tau}}w_\alpha(\vartheta)\,d\vartheta\right|^pw_\alpha(\tau)\, d\tau\right.\\ +
\left.\int_\alpha^{\pi}\left|\psi_{n+1}(e^{i\tau})\int_\alpha^{\pi+\delta} f(e^{i\vartheta})\frac{\psi_{n+1}^*(e^{i\alpha})}{\psi_{n+1}(e^{i\alpha})}\frac{\overline{\Lambda_{n+1}^{(\alpha)}(e^{i\vartheta)})}}{1-\overline{e^{i\vartheta}}e^{i\tau}}w_\alpha(\vartheta)\,d\vartheta\right|^pw_\alpha(\tau)\right).
\end{multline}
According to Lemma \ref{LemaAcotLambda}, we know that there exists $C>0$ such that
\[
\left|\frac{\psi_{n+1}^*(e^{i\alpha})}{\psi_{n+1}(e^{i\alpha})}\overline{\Lambda_{n+1}^{(\alpha)}(e^{i\vartheta)})}w_\alpha(\vartheta)\right|\le C,\quad \text{for all }\vartheta\in(\alpha,\pi+\delta),
\]
and from \eqref{eqBoundAkh} we know that
\[
|\psi_{n+1}(e^{i\tau})|\le C,\quad \text{for all }\tau \in (0,\pi).
\]
Therefore, by Lemma \ref{MucOnAnArc} there exists $C>0$ such that for all $f\in L^p(w_\alpha)$ we have
\begin{multline}
\label{DesPrimera}
\int_\alpha^{\pi}\left|\psi_{n+1}(e^{i\tau})\int_\alpha^{\pi+\delta} f(e^{i\vartheta})\frac{\frac{\psi_{n+1}^*(e^{i\alpha})}{\psi_{n+1}(e^{i\alpha})}\overline{\Lambda_{n+1}^{(\alpha)}(e^{i\vartheta)})}}{1-\overline{e^{i\vartheta}}e^{i\tau}}w_\alpha(\vartheta)\,d\vartheta\right|^pw_\alpha(\tau)\\ \le C\int_{\alpha}^{2\pi-\alpha}\left|f(e^{i\vartheta})\right|^pw_\alpha(\vartheta)\, d\vartheta.
\end{multline}

On the other hand, doing the change of variables $e^{i\tau}=h(e^{is})$, $e^{i\vartheta}=h(e^{it})$ and $s,t\in (0,\pi/2)$ and by the first equality in \eqref{MainTransfH},
\begin{multline*}
\int_\alpha^{\pi}\left|\Lambda_{n+1}^{(\alpha)}(e^{i\tau})\int_\alpha^{\pi+\delta} f(e^{i\vartheta})\frac{\overline{\psi_{n+1}^*(e^{i\vartheta)})}}{1-\overline{e^{i\vartheta}}e^{i\tau}}w_\alpha(\vartheta)\,d\vartheta\right|^pw_\alpha(\tau)\, d\tau \\ \le C\int_0^{\pi/2}\left|\frac{\Lambda_{n+1}^{(\alpha)}(h(e^{is}))}{\sin s}\int_0^{\pi/2+\widehat{\delta}} f(h(e^{it}))\overline{\psi_{n+1}^*(h(e^{it}))} (\cot\frac{t+s}{2}-\cot\frac{t-s}{s})\,dt\right|^pds,
\end{multline*}
where the arc in the unit circle from $1$ to $e^{i(\pi/2+\widehat{\delta})}$ has image by $h$  the arc from $e^{i\alpha}$ to $e^{i(\pi+\delta)}$. From \eqref{FormulaSin} and Lemma  \ref{LemaAcotLambda} we know that there exists a constant $C$ such that
\[
\left|\frac{\Lambda_{n+1}^{(\alpha)}(h(e^{is}))}{\sin s}\right|\le C, \quad \text{for all } s\in (0,\frac{\pi}{2}).
\]
Then using Riesz's result \eqref{bouHilbTrans2}, it follows the inequality
\begin{multline}
\label{DesSegunda}
\int_\alpha^{\pi}\left|\Lambda_{n+1}^{(\alpha)}(e^{i\tau})\int_\alpha^{\pi+\delta} f(e^{i\vartheta})\frac{\overline{\psi_{n+1}^*(e^{i\vartheta)})}}{1-\overline{e^{i\vartheta}}e^{i\tau}}w_\alpha(\vartheta)\,d\vartheta\right|^pw_\alpha(\tau)\, d\tau \\ \le C
\int_\alpha^{2\pi-\alpha}|f(e^{i\tau})|^pw_\alpha(\tau)\, d\tau.
\end{multline}
Plugging \eqref{DesPrimera} and \eqref{DesSegunda} into \eqref{RelaBasicaDiagonal}, we  finish the proof of \eqref{AcotaFourierSeriesPsi1}.
\end{proof}

\section{Additional results}\label{SectionAdd}

This section contains some results  about mean convergence of series in orthogonal  polynomials with respect to Akhiezer-Chebyshev weight multiplied by a function with nice properties. Assume that $k(\vartheta)\ge k>0$ for all $\vartheta\in(\alpha,2\pi-\alpha)$ and let $k(\vartheta)$ satisfy the Lipschitz condition in $\Delta_\alpha$, i.e.
\begin{equation}
\label{LipCond} 
|k(\vartheta_1)-k(\vartheta_2)|\le \lambda |\vartheta_1-\vartheta_2|,\quad \vartheta_1,\vartheta_2\in(\alpha,2\pi-\alpha) ,
\end{equation}
where $\lambda$ is a positive constant. We consider the measure
\begin{equation}
\label{AkhCheTypWei}
d\mu_{\alpha}(\vartheta)=k(\vartheta)w_\alpha(\vartheta)\, d\vartheta,\quad \vartheta\in(\alpha,2\pi-\alpha).
\end{equation}
 Observe that since the function $k(\vartheta)$ satisfies the Lipschitz condition, it is also bounded above in the interval $(\alpha,2\pi-\alpha)$.

To obtain our main results in this section we shall need some auxiliary results and definitions. The weighted Hilbert transform  is defined by
 \begin{equation*}
\mathcal{H}_2(f)(e^{i\tau})=\mathcal{H}_1(fw_\alpha)(e^{i\tau})
=\frac{1}{2\pi}\lim_ {\epsilon\to 0^+}\int_{(0,\pi)\setminus(s-\epsilon,s+\epsilon)}\frac{f(h(e^{it}))}{h(e^{it})-h(e^{is})}dt,
\end{equation*}
where $\tau\in(\alpha,2\pi-\alpha)$, $e^{i\tau}=h(e^{is})$, $e^{i\vartheta}=h(e^{it})$ and $s,t\in (0,\pi)$.

For the weighted Hilbert transform we have the following inequality.
\begin{lemma}\label{lemHilbert}  There exists a constant $C$ such that for all $f\in L^p(w_\alpha)$ we have
\begin{multline*}
\int_\alpha^{2\pi-\alpha}\left|\mathcal{H}_2(f)(e^{i\tau})\sqrt{\cos^2(\alpha/2)-\cos^2(\tau/2)}\right|^pw_\alpha(\tau)\, d\vartheta\\ \le C\int_\alpha^{2\pi-\alpha}|f(e^{i\tau})|^p\, w_\alpha(\tau)\, d\tau.
\end{multline*}
\end{lemma}

\begin{proof}  
Following the same steps of the proof of  Lemma \ref{MucOnAnArc}   until \eqref{eqHilbertTransf1} we get
\begin{multline*}
\int_\alpha^{2\pi-\alpha}|\mathcal{H}_2(f)(e^{i\tau})\sqrt{\cos^{2}(\alpha/2)-\cos^2(\tau/2)}|^pw_\alpha(\tau)\, d\tau\\=\int_\alpha^{2\pi-\alpha}\left|\int_\alpha^{2\pi-\alpha}\frac{f(e^{i\vartheta})\sqrt{\cos^{2}(\alpha/2)-\cos^2(\tau/2)}}{e^{i\vartheta}-e^{i\tau}}w_\alpha(\vartheta)\, d\vartheta\right|^pw_\alpha(\tau)\, d\tau\\ =\int_0^{\pi}\left|\int_0^{\pi}\frac{f(h(e^{it}))\sin s\cos(\alpha/2)(\sin\tau/2)}{h(e^{is})-h(e^{it})}\, dt\right|^p\, ds.
\end{multline*}
Then we obtain
\begin{multline}
\label{eqHilbertTransf}
\left(\int_\alpha^{2\pi-\alpha}|\mathcal{H}_2(f)(e^{i\tau})\sqrt{\cos^{2}(\alpha/2)-\cos^2(\tau/2)}|^pw_\alpha(\tau)\, d\tau\right)^{1/p}\\\le C \left(\int_0^\pi |\overset{\backsim}{f^*}(s)|^p\, ds\right)^{1/p}
\end{multline}
where $f^*(s)=f(h(e^{is}))$ and the constant $C>0$ is independent of $f$. Therefore, the lemma follows from \eqref{bouHilbTrans2}, \eqref{eqHilbertTransf} and Lemma~\ref{LemConfMap}~\emph{(iii)}.
\end{proof}

Finally, we need to state a type of Korus' lemma (\cite[p. 162]{Sze75}).

\begin{lemma}\label{LemAuxBoundPol} Let $\mu_\alpha$ be the measure given by \eqref{AkhCheTypWei} with $k(\vartheta)$ such that $k(\vartheta)\ge k>0$ for all $\vartheta\in(\alpha,2\pi-\alpha)$ and let $k(\vartheta)$ satisfy the Lipschitz condition in $\Delta_\alpha$ given by \eqref{LipCond}. Let $\{\varphi_{n}\}_{n=0}^\infty$ denote sequence of ortonormal polynomials with respect to $\mu_\alpha$. Then there exists a constant $C$ such that
\[
|\varphi_{n}(e^{i\vartheta})|\le C,\quad\text{for all } \vartheta\in (\alpha,2\pi-\alpha),\, \text{and all }n\in\N.
\]
\end{lemma}

\begin{proof} Let $\kappa_n$ denote the leading coefficient of $\varphi_n$. By the reproducing property of Christoffel kernel and Christoffel-Darboux formula we have
\begin{multline*}
\varphi_{n}(e^{i\tau})=\int_{\alpha}^{2\pi-\alpha} \varphi_{n}(e^{i\vartheta})\sum_{j=0}^n\psi_j(e^{i\tau})\overline{\psi_j(e^{i\vartheta})}w_\alpha(\vartheta)\, d\vartheta\\=\frac{\kappa_n}{\alpha_n}\psi_{n}(e^{i\tau})+\int_{\alpha}^{2\pi-\alpha}\varphi_{n}(e^{i\vartheta})\sum_{j=0}^{n-1}\psi_j(e^{i\tau})\overline{\psi_j(e^{i\vartheta})}w_\alpha(\vartheta)\left(1-\frac{k(\vartheta)}{k(\tau)}\right)\, d\vartheta\\=\frac{\kappa_n}{\alpha_n}\psi_{n}(e^{i\tau})\\+\frac{1}{k(\tau)}\int_{\alpha}^{2\pi-\alpha} \varphi_{n}(e^{i\vartheta})\left(\overline{\psi_n^*(e^{i\vartheta})}\psi_n^*(e^{i\tau})-\overline{\psi_n(e^{i\vartheta})}\psi_n(e^{i\tau})\right)\frac{k(\tau)-k(\vartheta)}{1-e^{i(\tau-\vartheta)}}w_\alpha(\vartheta)\, d\vartheta.
\end{multline*}
Since\footnote{As usual, the notation $a(x)\sim b(x)$ for $x$ in an interval $I$ means that there exist positive constants $C_1,C_2$ independent of $x$ such that $C_1\le \frac{a(x)}{b(x)}\le C_2$ for all $x\in I$.}
\[
|1-e^{i(\tau-\vartheta)}|= 2|\sin\frac{\tau-\vartheta}{2}|\sim |\tau-\vartheta|,\quad \alpha<\tau,\vartheta<2\pi-\alpha,
\]
by \eqref{LipCond}, there exists a constant $C>0$ such that
\[
\left|\frac{k(\tau)-k(\vartheta)}{1-e^{i(\tau-\vartheta)}}\right|\le C.
\]
By Cauchy-Schwarz's inequality
\begin{multline*}
\frac{\kappa_n}{\alpha_n}=\int_{\alpha}^{2\pi-\alpha}\varphi_{n}(e^{i\vartheta})\overline{\psi_n(e^{i\vartheta})}w_\alpha(\vartheta)\, d\vartheta \\ \le \left(\int_{\alpha}^{2\pi-\alpha} |\varphi_{n}(e^{i\vartheta})|^2\, \frac{d\mu_\alpha(\vartheta) }{k(\vartheta)}\right)^{1/2}\left(\int_{\alpha}^{2\pi-\alpha} |\psi_n(e^{i\vartheta})|^2w_\alpha(\vartheta)\, d\vartheta \right)^{1/2}\le k^{-1/2}.
\end{multline*}
Taking into account Lemma \ref{lemAux} we known that $\{\psi_n(e^{i\vartheta})\}_{n=0}^{\infty}$ is uniformly bounded for $\vartheta\in (\alpha,2\pi-\alpha)$, therefore the above inequalities prove the lemma.
\end{proof}

Let $S_n(f,\cdot)$ denote the $n-$th partial  Fourier sums in terms of $\{\varphi_j\}_{j=0}^\infty$ for the function $f\in L^p(\mu_\alpha)$.

\begin{theorem}\label{TeoConvMedSerieFourier} Let $\mu_\alpha$ be the measure given by \eqref{AkhCheTypWei} with $k(\vartheta)$ such that $k(\vartheta)\ge k>0$ for all $\vartheta\in(\alpha,2\pi-\alpha)$ and let $k(\vartheta)$ satisfy the Lipschitz condition in $\Delta_\alpha$ given by \eqref{LipCond}. Then for all $f\in L^p(\mu_\alpha)$, $1<p<\infty$, we have
\[
\lim_{n\to\infty}\int_{\alpha}^{2\pi-\alpha}|(f(e^{i\tau})-S_n(f,e^{i\tau}))\sqrt{\cos^2(\alpha/2)-\cos^2(\tau/2)}|^p\, d\mu_\alpha(\tau)=0.
\]
\end{theorem}

\begin{proof}

According to Christoffel-Darboux formula  
\begin{multline}\label{ChrDarFor}
S_n(f,z)=   \frac{\varphi_{n+1}^*(z)}{2\pi}\int_{\alpha}^{2\pi-\alpha} \frac{f(e^{i\vartheta})\overline{\varphi_{n+1}^*(e^{i\vartheta})}}{1-\overline{e^{i\vartheta}}z}\,d\mu_{\alpha}(\vartheta)\\-\frac{\varphi_{n+1}(z)}{2\pi}\int_{\alpha}^{2\pi-\alpha}\frac{f(e^{i\vartheta})\overline{\varphi_{n+1}(e^{i\vartheta})}}{1-\overline{e^{i\vartheta}}z}\,d\mu_{\alpha}(\vartheta).
\end{multline}
Because of Lemma \ref{LemAuxBoundPol} we know that the Akhiezer-Chebyshev type polynomials $\{\varphi_n\}_{n=0}^{\infty}$ is a uniformly bounded sequence  on the arc $\Delta_\alpha$.  Combining Lemma \ref{lemHilbert}, the hypothesis on the function $k$, and \eqref{ChrDarFor}, for all $p>1$ there exists a constant $C=C(p)>0$ such that for all $f\in L^p(w_\alpha)$ we have
\begin{equation}
\label{AcotSumFourier}
\int_\alpha^{2\pi-\alpha}|S_n(f,e^{i\tau})\sqrt{\cos^{2}(\alpha/2)-\cos^2(\tau/2)}|^p\, d\mu_\alpha(t)\le C \int_0^\pi|f(e^{i\tau})|^p\, d\mu_\alpha(\tau),
\end{equation}
where $C$ is a constant independent of $f\in L^p(\mu_\alpha)$. From the Szeg\H{o}-Kolmogorov-Krein theorem, the algebraic  polynomials are dense in $L^p(\mu_\alpha)$. Thus, there exists a sequence of algebraic polynomials $\{p_n\}$ with $\text{deg}(p_n)\le n$ such that
\[
\lim_{n\to\infty}\int_{\alpha}^{2\pi-\alpha}|f(e^{i\tau})-p_n(e ^{i\tau})|^p\, d\mu_\alpha(\tau)=0.
\]
By \eqref{AcotSumFourier}, we have
\begin{multline*}
\left(\int_\alpha^{2\pi-\alpha}|(S_n(f,e^{i\tau})-f(e^{i\tau}))\sqrt{\cos^{2}(\alpha/2)-\cos^2(\tau/2)}|^p\, d\mu_\alpha(\vartheta)\right)^{1/p}
\\ \le 
\left(\int_\alpha^{2\pi-\alpha}|(S_n(f-p_n,e^{i\tau}))\sqrt{\cos^{2}(\alpha/2)-\cos^2(\tau/2)}|^p\, d\mu_\alpha(\tau)\right)^{1/p}\\\quad+\left(\int_\alpha^{2\pi-\alpha}|(p_n(e^{i\tau})-f(e^{i\tau}))\sqrt{\cos^{2}(\alpha/2)-\cos^2(\tau/2)}|^p\, d\mu_\alpha(\tau)\right)^{1/p}
\\ \le C\left(\int_{\alpha}^{2\pi-\alpha}|f(e^{i\tau})-p_n(e ^{i\tau})|^p\, d\mu_\alpha(\tau)\right)^{1/p},
\end{multline*}
and the conclusion of the theorem follows.
\end{proof}

\begin{remark} Since the operator which maps $f\in L^\infty(\mu_ \alpha)$ into its $n-$th partial Fourier series is a projection operator, according to Losinski-Kharshilarze-Nikolaev's theorem \cite[Appendix 3]{Nat65}, this operator is not bounded from $L^\infty(\mu_ \alpha)$ to $L^\infty(\mu_ \alpha)$. By duality, it is also not   bounded from $L^1(\mu_\alpha)$ to itself.

Of course, Theorems \ref{TeoConvMedSerieFourier} is not sharp for $p=2$.
An improvement of Theorem \ref{TeoConvMedSerieFourier} would be obtained if an inequality like in Lemma \ref{LemaAcotLambda}  is proved for para-orthogonal polynomials associated to the measure $\mu_\alpha$.

\end{remark}

\vspace{0.5cm}

\textbf{Acknowledgments.} The authors would like to thank the referees for helping us to improve the presentation of this paper and draw our attention to the reference \cite{SteWai78}.

\bibliographystyle{acm}


\end{document}